\documentclass[oneside,english]{amsart}
\usepackage[T1]{fontenc}
\usepackage[latin9]{inputenc}
\usepackage{geometry}
\usepackage{amsthm}
\usepackage{amstext}
\usepackage{amssymb}
\usepackage{esint}

\newtheorem{theorem}{Theorem}

\newtheorem{rem}[theorem]{Remark}

\newtheorem{thm}[theorem]{Theorem}

\makeatletter

\makeatletter
\numberwithin{equation}{section}

\usepackage{babel}

\author{R. M. \L OCHOWSKI} 

\makeatother

\usepackage{babel}
\begin{document}

\title{On pathwise uniform approximation \\
 of processes with càdlàg trajectories \\
 by processes with finite total variation}

\maketitle
\textbf{Abstract.} For any real-valued stochastic process $X=\left(X_{t}\right)_{t\geq0}$
with càdlàg paths we define non-empty family of processes, which have
finite total variation, have jumps of the same order as the process
$X$ and uniformly approximate its paths$.$\ This allows to decompose
any real-valued stochastic process with càdlàg paths and infinite
total variation into a sum of uniformly close, finite variation process
and an adapted process, with arbitrary small amplitude but infinite
total variation. Another application of the defined class is the definition
of the stochastic integral with respect to the process $X$ as a limit
of pathwise Lebesgue-Stieltjes integrals. This construction leads
to the stochastic integral with some correction term.

\section{Introduction }

Let $X=\left(X_{t}\right)_{t\geq0}$ be a real-valued stochastic process
with càdlàg paths and let $0\leq a<b.$ The total variation of the
process $X$ on the interval $\left[a;b\right]$ is defined with the
following formula 
\[
TV\left(X,\left[a;b\right]\right)=\sup_{n}\sup_{a\leq t_{0}<t_{1}<...<t_{n}\leq b}\sum_{i=1}^{n}\left|X_{t_{i}}-X_{t_{i-1}}\right|.
\]
 Unfortunately, many of the most important families of stochastic
processes are characterized with the ''wild'' behaviour, demonstrated
by their infinite total variation. This fact arguably caused the need
of the development of the general theory of stochastic integral. The
main idea allowing to overcome the problematic infinite total variation
and define stochastic integral with respect to a semimartingale utilizes
the fact that the quadratic variation of the semimartingale is still
finite. The similar idea may be applied when $p-$variation of the
integrator is finite for some $p>1$. This approach utilizes Love-Young
inequality and may be used e.g. to define stochastic integral with
respect to fractional Brownian motion (cf. \cite{Kubilius:2008}).
Further developments, where Hölder continuity plays crucial role,
led to the rough paths theory developed by T. Lyons and his co-workers
(cf. \cite{Friz:2010fk}); some other generalization introduces Orlicz
norms and may be found in the recent book by Dudley and Norvaiša \cite[Chapt. 3]{DudleyNorvaisa:2011}).
The approach used in this article is somewhat different. It is similar
to the old approach of Wong and Zakai \cite{WongZakai:1965a} and
is based on the simple observation that in the neighborhood (in sup
norm) of every càdlàg function defined on some compact interval one
easily finds another function with finite total variation. Thus, for
every $c>0,$ the process $X$ may be decomposed as the sum 
\[
X=X^{c}+\left(X-X^{c}\right)
\]
 where $X^{c}$ is a {}``nice'' process with finite total variation
and the difference $X-X^{c}$ is a process with small amplitude (no
greater than $c$) but possibly {}``wild'' behaviour with infinite
total variation. More precisely, let $F$ be some fixed, right continuous
filtration such that $X$ is adapted to $F.$ Now, for every $c>0$
we introduce (non-empty, as it will be shown in the sequel) family
$\mathcal{X}^{c}$ of processes with càdlàg paths, satisfying the
following conditions. If $X^{c}\in\mathcal{X}^{c}$ then 
\begin{enumerate}
\item the process $X^{c}$ has locally finite total variation; 
\item $X^{c}$ has càdlàg paths; 
\item for every $t\geq0,$ $\left|X_{t}-X_{t}^{c}\right|\leq c;$ 
\item for every $T\geq0$ there exists such $K_{T}<+\infty$ that for every
$t\in\left[0;T\right],$ $\left|\Delta X_{t}^{c}\right|\leq K_{T}\left|\Delta X_{t}\right|$$;$ 
\item the process $X^{c}$ is adapted to the filtration $F.$ 
\end{enumerate}
We will prove that if processes $X$ and $Y$ are càdlàg semimartingales
on a filtered probability space $\left(\Omega,{\mathcal{F}},\mathbb{P},F\right),$
with a probability measure $\mathbb{P},$ such that usual hypotheses
hold (cf. \cite[Sect. 1.1]{Protter:2004}), then the sequence of pathwise
Lebesgue-Stieltjes integrals 
\[
\int_{0}^{T}Y_{-}\mathrm{d}X^{c},\, c>0,
\]
 with $X^{c}\in{\mathcal{X}}^{c}$, tends uniformly in probability
$\mathbb{P}$ on compacts to $\int_{0}^{T}Y_{-}\mathrm{d}X+\left[X^{cont},Y^{cont}\right]_{T}$;
$\int_{0}^{T}Y_{-}\mathrm{d}X$ denotes here the (semimartingale)
stochastic integral and $X^{cont}$ and $Y^{cont}$ denote continuous
parts of $X$ and $Y$ respectively. Moreover, for any square summable
sequence $\left(c\left(n\right)\right)_{n\geq1}$ we get $\mathbb{P}$ a.s. and uniform on compacts convergence of the sequence $\int_{0}^{T}Y_{-}\mathrm{d}X^{c\left(n\right)},n=1,2,...$
(cf. Theorem \ref{thm:leb_stieltjes_stoch-1}). We shall stress here
that for every $c>0$ and each pair of càdlàg paths $\left(X\left(\omega\right),Y\left(\omega\right)\right),\omega\in\Omega,$
the value of $\int_{0}^{T}Y_{-}\left(\omega\right)\mathrm{d}X^{c}\left(\omega\right)$
(and thus the limit, if it exists) is independent of the probability
measure $\mathbb{P}$. Thus we obtain a result in the spirit of Bichtelier,
see \cite{Karandikar:1995}, and recent result of Nutz \cite{Nutz:2011},
where operations leading to the stochastic integral, independent of
probability measures and filtrations are considered. Our approach
seems to be simpler and more natural, however we need to impose a
stronger condition on the integrand - that it is also a semimartingale.

Further, for $p\geq1$ we will also investigate the behaviour of $p-$variation
of the processes $X^{c}$ and $X-X^{c}$ as $c\downarrow0.$ E.g.,
for $X$ being a semimartingale, the properties (1)-(5) allow to fully
determine the (almost sure) limits 
\[
\lim_{c\downarrow0}v_{p}\left(X^{c};\left[0;T\right]\right)
\]
 in terms of predictable characteristics of $X$ (for definition of
predictable characteristics see \cite{JacodSh:2003}). $v_{p}\left(Z;\left[0;T\right]\right)$
denotes here $p-$variation of a process $Z$ on an interval $\left[0;T\right],$
defined in few different ways (see Section 4 for details). The limits
$\lim_{c\downarrow0}v_{p}\left(X^{c};\left[0;T\right]\right)$ will
coincide with the pathwise limit of equidistant $p-$variations of
the process $X^{c},$ defined as 
\begin{equation}
v_{p}^{\left(0\right)}\left(X^{c};\left[0;T\right]\right):=\limsup_{n\rightarrow\infty}\sum_{i=1}^{n}\left|X_{iT/n}^{c}-X_{\left(i-1\right)T/n}^{c}\right|.\label{eq:realized}
\end{equation}
 The investigation of limits of equidistant $p-$variation of stochastic
processes, possibly perturbed with some noise, as the mesh of the
partitions goes to $0$, may be of practical interest. E.g., functional
limits of equidistant $p-$variation of $\alpha-$ stable processes,
perturbed with some noise, were investigated and used in \cite{HeinImkeller:2009}
to model paleoclimatic temperature time series taken from the Greenland
ice core. In \cite{AitSahaliaJacod:2009} limits of equidistant $p-$variation
were investigated and used for testing whether jumps are present in
asset returns or other, discretely sampled processes.

Let us shortly comment on the organization of the paper. In the next
section we prove, for any $c>0,$ the existence of non-empty family
of processes ${\mathcal{X}}^{c}$. In the third section we deal with
the limit of pathwise, Lebesgue-Stieltjes integrals $\int_{0}^{T}Y_{-}\mathrm{d}X^{c}$
as $c\downarrow0$ and in the last section we deal, for $p\geq1,$
with $p-$ variations of the processes $X^{c}$and $X-X^{c}.$

\textbf{Acknowledgments.} The author would like to thank Prof. Krzysztof
Burdzy for encouraging him to submit this paper by saying that the
problems considered are interesting and to thank Dr. Alexander Cox
for pointing out to him the results of \cite{Nutz:2011}.

\section{Existence of the sequence $\left(X^{c}\right)_{c>0}$}

In this section we will prove that for every $c>0$ the family of
processes ${\mathcal{X}}^{c},$ satisfying the conditions (1)-(5)
of Section 1 is non-empty. For given $c>0$ we will simply construct
a process $X^{c}$ satisfying all these conditions. Our construction
is neither unique nor optimal (in the sense that it does not produce
a process satisfying (1)-(5) with the smallest total variation possible),
but it seems to be the simplest one. We start with few definitions.

For fixed $c>0$ we define two stopping times 
\begin{gather*}
T_{u}^{2c}X=\inf\left\{ s\geq0:\mbox{ }\sup_{t\in\left[0;s\right]}X_{t}-X_{0}>c\right\} ,\\
T_{d}^{2c}X=\inf\left\{ s\geq0:\mbox{ }X_{0}-\inf_{t\in\left[0;s\right]}X_{t}>c\right\} .
\end{gather*}

Assume that $T_{d}^{2c}X\geq T_{u}^{2c}X,$ i.e. the first upward
jump of the process $X$ from $X_{0}$ of size $c$ appears before
the first downward jump of the same size $c$ or both times are infinite
(there is no upward or downward jump of size $c$). Note that in the
case $T_{d}^{2c}X<T_{u}^{2c}X$ we may simply consider the process
$-X.$ Now we define sequences $\left(T_{d,k}^{2c}\right)_{k=1}^{\infty},\mbox{ }\left(T_{u,k}^{2c}\right)_{k=1}^{\infty}$
in the following way: $T_{u,0}^{2c}=T_{u}^{2c}X$ and for $k=0,1,2,...$
\begin{gather*}
T_{d,k}^{2c}=\left\{ \begin{array}{lr}
\inf\left\{ s\geq T_{u,k}^{2c}:\sup_{t\in\left[T_{u,k}^{2c};s\right]}X_{t}-X_{s}>2c\right\}  & \text{ if }T_{u,k}^{2c}<+\infty,\\
+\infty & \mbox{otherwise},
\end{array}\right.\\
T_{u,k+1}^{2c}=\left\{ \begin{array}{lr}
\inf\left\{ s\geq T_{d,k}^{2c}:X_{s}-\inf_{t\in\left[T_{d,k}^{2c};s\right]}X_{t}>2c\right\}  & \text{ if }T_{d,k}^{2c}<+\infty,\\
+\infty & \mbox{otherwise}.
\end{array}\right.
\end{gather*}
 \begin{rem} \label{finK} Note that for any $s>0$ there exists
such $K<\infty$ that $T_{u,K}^{2c}>s$ or $T_{d,K}^{2c}>s.$ Otherwise
we would obtain two infinite sequences $\left(s_{k}\right)_{k=1}^{\infty},\left(S_{k}\right)_{k=1}^{\infty}$
such that $0\leq s\left(1\right)<S\left(1\right)<s\left(2\right)<S\left(2\right)<...\leq s$
and $X_{S\left(k\right)}-X_{s\left(k\right)}\geq c.$ But this is
a contradiction since $X$ is a càdlàg process and for any sequence
such that $0\leq s\left(1\right)<S\left(1\right)<s\left(2\right)<S\left(2\right)<...\leq s$
sequences $\left(X_{S\left(k\right)}\right)_{k=1}^{\infty},\left(X_{s\left(k\right)}\right)_{k=1}^{\infty}$
have a common limit. \end{rem}

Now we define, for the given process $X,$ the process $X^{c}$ with
the formulas 
\begin{equation}
X_{s}^{c}=\left\{ \begin{array}{lr}
X_{0} & \text{ if }s\in\left[0;T_{u,0}^{2c}\right);\\
\sup_{t\in\left[T_{u,k}^{2c};s\right]}X_{t}-c & \text{ if }s\in\left[T_{u,k}^{2c};T_{d,k}^{2c}\right),k=0,1,2,...;\\
\inf_{t\in\left[T_{d,k}^{2c};s\right]}X_{t}+c & \text{ if }s\in\left[T_{d,k}^{2c};T_{u,k+1}^{2c}\right),k=0,1,2,....
\end{array}\right.\label{eq:defXc}
\end{equation}

\begin{rem} Note that due to Remark \ref{finK}, $s$ belongs to
one of the intervals $\left[0;T_{u,0}^{2c}\right),\left[T_{u,k}^{2c};T_{d,k}^{2c}\right)$
or $\left[T_{d,k}^{2c};T_{u,k+1}^{c}\right)$ for some $k=0,1,2,...$
and the process $X_{s}^{c}$ is defined for every $s\geq0.$ \end{rem}

Now we are to prove that $X^{c}$ satisfies conditions (1)-(5).

\begin{proof}

(1) The process $X^{c}$ has finite total on compact intervals, since
it is monotonic on intervals of the form $\left[T_{u,k}^{2c};T_{d,k}^{2c}\right),$
$\left[T_{d,k}^{2c};T_{u,k+1}^{c}\right)$ which sum up to the whole
half-line $\left[0;+\infty\right).$

(2) From formula (\ref{eq:defXc}) it follows that $X^{c}$ is also
càdlàg.

(3) In order to prove condition (3) we consider 3 possibilities. 
\begin{itemize}
\item $s\in\left[0;T_{u,0}^{2c}\right).$ In this case, since $0\leq s<T_{u}^{2c}X\leq T_{d}^{2c}X,$
by definition of $T_{u}^{2c}X$ and $T_{d}^{2c}X,$ 
\[
X_{s}-X_{s}^{c}=X_{s}-X_{0}\in\left[-c;c\right].
\]
 
\item $s\in\left[T_{u,k}^{2c};T_{d,k}^{2c}\right),$ for some $k=0,1,2,...$
In this case, by definition of $T_{d,k}^{2c}$, $\sup_{t\in\left[T_{u,k}^{2c};s\right]}X_{t}-X_{s}$
belongs to the interval $\left[0;2c\right],$ hence 
\[
X_{s}-X_{s}^{c}=X_{s}-\sup_{t\in\left[T_{u,k}^{2c};s\right]}X_{t}+c\in\left[-c;c\right].
\]
 
\item $s\in\left[T_{d,k}^{2c};T_{u,k+1}^{2c}\right)$ for some $k=0,1,2,...$
In this case $X_{s}-\inf_{t\in\left[T_{d,k}^{2c};s\right]}X_{t}$
belongs to the interval $\left[0;2c\right],$ hence 
\[
X_{s}-X_{s}^{c}=X_{s}-\inf_{t\in\left[T_{d,k}^{2c};s\right]}X_{t}-c\in\left[-c;c\right].
\]
 
\end{itemize}
(4) We will prove stronger fact than (4), namely that for every $s>0,$
\begin{align}
\left|\Delta X_{s}^{c}\right| & \leq\left|\Delta X_{s}\right|.\label{eq:jumps}
\end{align}
Indeed, from formula (\ref{eq:defXc}) it follows that for any $s\notin\left\{ T_{u,k}^{2c};T_{d,k}^{2c}\right\} ,$
(\ref{eq:jumps}) holds, hence let us assume that $s\in\left\{ T_{u,k}^{2c};T_{d,k}^{2c}\right\} .$
We consider several possibilities. If $s=T_{u,0}^{2c}$ then, by the
definition of $T_{u,0}^{2c},$ 
\[
X_{s}^{c}-X_{s-}^{c}=X_{s}-c-X_{0}\mbox{ \ensuremath{\geq}0 and }X_{s}^{c}-X_{s-}^{c}=X_{s}-X_{0}-c\leq X_{s}-X_{s-}.
\]
 If $s=T_{u,k}^{2c},k=1,2,...,$ then, by the definition of $T_{u,k}^{2c},$
\[
X_{s}^{c}-X_{s-}^{c}=X_{s}-c-\left(\inf_{t\in\left[T_{d,k-1}^{2c};s\right]}X_{t}+c\right)\mbox{ =\ensuremath{X_{s}-\inf_{t\in\left[T_{d,k-1}^{2c};s\right]}X_{t}-2c\geq0}}
\]
 and, on the other hand, 
\[
X_{s}^{c}-X_{s-}^{c}\mbox{ =\ensuremath{X_{s}-\inf_{t\in\left[T_{d,k-1}^{2c};s\right]}X_{t}-2c\leq X_{s}-X_{s-}.}}
\]
 Similar arguments may be applied for $s=T_{d,k}^{2c},k=0,1,....$

(5) The process $X^{c}$ is adapted to the filtration $F$ since it
is adapted to any right continuous filtration containing the natural
filtration of the process $X$.

\end{proof} 

\begin{rem} It is possible to define the process $X^{c}$ in many
different ways. For example, defining 
\[
X^{c}=X_{0}+UTV^{c}\left(X,\cdot\right)-DTV^{c}\left(X,\cdot\right)
\]
 we obtain a process satisfying all conditions (1)-(5) and having
(on the intervals of the form $\left[0;T\right],\mbox{ }T>0$) the
smallest possible total variation among all processes, increments
of which differ from the increments of the process $X$ by no more
than $c.$ $UTV^{c}\left(X,\cdot\right)\mbox{ and }DTV^{c}$ denote
here upward and downward truncated variation processes, defined as
\begin{eqnarray*}
UTV^{c}\left(X,t\right) & = & \sup_{n}\sup_{0\leq t_{1}<t_{2}<...<t_{n}\leq t}\sum_{i=1}^{n}\max\left\{ X_{t_{i}}-X_{t_{i-1}}-c,0\right\} ,\\
DTV^{c}\left(X,t\right) & = & \sup_{n}\sup_{0\leq t_{1}<t_{2}<...<t_{n}\leq t}\sum_{i=1}^{n}\max\left\{ X_{t_{i-1}}-X_{t_{i}}-c,0\right\} .
\end{eqnarray*}
 Moreover, for any $T>0$ we have 
\begin{eqnarray*}
TV\left(X^{c};\left[0;T\right]\right) & = & UTV^{c}\left(X,T\right)+DTV^{c}\left(X,T\right)\\
 & = & \sup_{n}\sup_{0\leq t_{1}<t_{2}<...<t_{n}\leq t}\sum_{i=1}^{n}\max\left\{ \left|X_{t_{i}}-X_{t_{i-1}}\right|-c,0\right\} =:TV^{c}\left(X;T\right).
\end{eqnarray*}
 For more on truncated variation, upward truncated variation and downward
truncated variation see e.g. \cite{L2011c} or \cite{LM2011}. \end{rem}

\section{Pathwise Lebesgue-Stieltjes integration with respect to the processes
$X^{c}$}

Let us now consider a measurable space $\left(\Omega,{\mathcal{F}}\right)$
equipped with a right-continuous filtration $F$ and two processes
$X$ and $Y$ with càdlàg paths, adapted to $F.$ For $T>0$ and for
a sequence of processes $\left(X^{c}\right)_{c>0}$ with $X^{c}\in{\mathcal{X}}^{c}$
let us consider the sequence 
\begin{equation}
\int_{0}^{T}Y_{-}\mathrm{d}X^{c}.\label{eq:int_c}
\end{equation}
 The integral in (\ref{eq:int_c}) is understood in the pathwise,
Lebesgue-Stieltjes sense (recall that for any $c>0,$ $X^{c}$ has
bounded variation). We have 
\begin{thm} \label{thm:leb_stieltjes_stoch}Assume
that $\mathbb{P}$ is a probability measure on $\left(\Omega,{\mathcal{F}}\right)$
such that $X$ and $Y$ are semimartingales with respect to this measure
and filtration $F,$ which is complete under $\mathbb{P}$, then 
\[
\int_{0}^{T}Y_{-}\mathrm{d}X^{c}\rightarrow^{ucp\mathbb{P}}\int_{0}^{T}Y_{-}\mathrm{d}X+\left[X^{cont},Y^{cont}\right]_{T}\,\mbox{as }c\downarrow0,
\]
where {}``$\rightarrow^{ucp\mathbb{P}}$'' denotes uniform convergence
on compacts in probability $\mathbb{P}$ and $\left[X^{cont},Y^{cont}\right]_{T}$
denotes quadratic covariation of continuous parts $X^{cont},Y^{cont}$
of $X$ and $Y$ respectively.\end{thm} 

\begin{proof} Fixing $c>0$ and using integration by parts formula
(cf. \cite[formula (1), page 519]{Kallenberg:2002}) we get 
\[
Y_{T}X_{T}^{c}-Y_{0}X_{0}^{c}=\int_{0}^{T}Y_{t-}\mathrm{d}X_{t}^{c}+\int_{0}^{T}X_{t-}^{c}\mathrm{d}Y_{t}+\left[Y,X^{c}\right]_{T}
\]
 (the above equality and subsequent equalities in the proof hold $\mathbb{P}$
a.s.). By the uniform convergence, $X_{t}^{c}\rightrightarrows X_{t}$
as $c\downarrow0$ (note that the bound $\left|X^{c}\right|\leq\left|X\right|+c$
and a.s. pointwise convergence $X_{t}^{c}\rightarrow X_{t}$ as $c\downarrow0$
are sufficient) we get 
\[
\int_{0}^{T}X_{t-}^{c}\mathrm{d}Y_{t}\rightarrow^{ucp\mathbb{P}}\int_{0}^{T}X_{t-}\mathrm{d}Y_{t}.
\]
 Since $X^{c}$ has locally finite variation, we have (cf. \cite[Theorem 26.6 (viii)]{Kallenberg:2002}),
\[
\left[Y,X^{c}\right]_{T}=\sum_{0<s\leq T}\Delta Y_{s}\Delta X_{s}^{c}.
\]
 We calculate the (pathwise) limit 
\[
\lim_{c\downarrow0}\left[Y,X^{c}\right]_{T}=\lim_{c\downarrow0}\sum_{0<s\leq T}\Delta Y_{s}\Delta X_{s}^{c}=\sum_{0<s\leq T}\Delta Y_{s}\Delta X_{s}
\]
 (notice that for any $0\leq s\leq T,$ $\left|\Delta X_{s}^{c}\right|\leq K_{T}\left|\Delta X_{s}\right|,$
thus the above sum is convergent by dominated convergence) and finally
obtain 
\begin{eqnarray}
\int_{0}^{T}Y_{t-}\mathrm{d}X_{t}^{c} & = & \left\{ Y_{T}X_{T}^{c}-Y_{0}X_{0}^{c}-\int_{0}^{T}X_{t-}^{c}\mathrm{d}Y_{t}-\left[Y,X^{c}\right]_{T}\right\} \nonumber \\
 & \rightarrow^{ucp\mathbb{P}} & Y_{T}X_{T}-Y_{0}X_{0}-\int_{0}^{T}X_{t-}\mathrm{d}Y_{t}-\sum_{0<s\leq T}\Delta Y_{s}\Delta X_{s}\mbox{ as }\ensuremath{c\downarrow0.}\label{eq:lim_int_part}
\end{eqnarray}
On the other hand, again by the integration by parts formula, we obtain
\begin{equation}
\int_{0}^{T}X_{t-}\mathrm{d}Y_{t}=Y_{T}X_{T}-Y_{0}X_{0}-\int_{0}^{T}Y_{t-}\mathrm{d}X_{t}-\left[Y,X\right]_{T}.\label{eq:int_part2}
\end{equation}
 Finally, comparing (\ref{eq:lim_int_part}) and (\ref{eq:int_part2}),
and using \cite[Corollary 26.15]{Kallenberg:2002}, we obtain 
\begin{eqnarray*}
\int_{0}^{T}Y_{t-}\mathrm{d}X_{t}^{c} & \rightarrow^{ucp\mathbb{P}} & \int_{0}^{T}Y_{t-}\mathrm{d}X_{t}+\left[Y,X\right]_{T}-\sum_{0<s\leq T}\Delta Y_{s}\Delta X_{s}\mbox{ \mbox{ as }\ensuremath{c\downarrow0}}\\
 & = & \int_{0}^{T}Y_{t-}\mathrm{d}X_{t}+\left[X^{cont},Y^{cont}\right]_{T}.
\end{eqnarray*}
 \end{proof} 

\begin{rem} Assuming the existence of Mokobodzki's medial limits
(cf. \cite{Meyer:1973}), which one can not prove under standard Zermelo\textendash{}Fraenkel
set theory with the axiom of choice, Theorem \ref{thm:leb_stieltjes_stoch}
may be used to construct a universal process which coincides with
stochastic integral, with the correction term $\left[X^{cont},Y^{cont}\right]_{T},$
for a family of probability measures simultaneously. More precisely,
we consider a family of probability measures ${\mathcal{P}}$ on $\left(\Omega,{\mathcal{F}}\right)$
such that for each $\mathbb{P}\in{\mathcal{P}}$ the filtration $F$
is complete under $\mathbb{P}$ and $X$ and $Y$ are semimartingales
on the filtered probability space $\left(\Omega,{\mathcal{F}},F,\mathbb{P}\right).$
Considering any sequence $\int_{0}^{T}Y_{-}\mathrm{d}X^{c\left(n\right)},\mbox{ \ensuremath{n=1,2,...},}$
with $c\left(n\right)\downarrow0,$ and using Theorem \ref{thm:leb_stieltjes_stoch}
and \cite[Lemma 2.5]{Nutz:2011} we obtain that there exists a universal,
$F$ adapted càdlàg process $I\left(Y,\mathrm{d}X,\left[0;\cdot\right]\right),$
such that for all $\mathbb{P}\in{\mathcal{P}}$ and $T>0$ 
\[
I\left(Y,\mathrm{d}X,\left[0;T\right]\right)=\int_{0}^{T}Y_{-}\mathrm{d}X+\left[X^{cont},Y^{cont}\right]_{T}\,\mathbb{P}\mbox{ a.s.}
\]

\end{rem}

Note that to prove Theorem \ref{thm:leb_stieltjes_stoch} we did not
need the pathwise uniform convergence of the processes $X^{c}$ to
the process $X;$ we might simply use local boundedness and a.s. pointwise
convergence $X_{t}^{c}\rightarrow X_{t}$ as $c\downarrow0.$ Using
the pathwise uniform convergence of the sequence $\left(X^{c}\right)_{c>0}$
we are able to prove a bit stronger result. We have 

\begin{thm} \label{thm:leb_stieltjes_stoch-1}
Assume that $\mathbb{P}$ is a probability measure on $\left(\Omega,{\mathcal{F}}\right)$
such that $X$ and $Y$ are semimartingales with respect to this measure
and filtration $F,$ which is complete under $\mathbb{P},$ then 
for any $T>0$ and any sequence $\left(c\left(n\right)\right)_{n\geq1}$
such that $\sum_{n=1}^{\infty}c\left(n\right)^{2}<+\infty$ we have
\[
\lim_{n\rightarrow+\infty}\sup_{0\leq t\leq T}\left|\int_{0}^{t}Y_{-}\mathrm{d}X^{c\left(n\right)}-\int_{0}^{t}Y_{-}\mathrm{d}X-\left[X^{cont},Y^{cont}\right]_{t}\right|=0 \mbox{ }\mathbb{P} \mbox{ a.s.}
\]

\end{thm} 

\begin{proof} Using integration by parts formula and the inequality
$\left|X^{c}-X\right|\leq c,$ we estimate 
\begin{eqnarray*}
\lefteqn{ \left|\int_{0}^{t}Y_{-}\mathrm{d}X^{c}-\int_{0}^{t}Y_{-}\mathrm{d}X-\left[X^{cont},Y^{cont}\right]_{t}\right| }
\\
& = & \left|Y_{t}X_{t}^{c}-Y_{0}X_{0}^{c}-\sum_{0<s\leq t}\Delta Y_{s}\Delta X_{s}^{c}-\int_{0}^{t}X_{-}^{c}\mathrm{d}Y-\left(Y_{t}X_{t}-Y_{0}X_{0}-\sum_{0<s\leq t}\Delta Y_{s}\Delta X_{s}-\int_{0}^{t}X_{-}^{c}\mathrm{d}Y\right)\right| 
\\ 
& = &
\left|Y_{t}\left(X_{t}^{c}-X_{t}\right)-Y_{0}\left(X_{0}^{c}-X_{0}\right)-\sum_{0<s\leq t}\Delta Y_{s} \Delta\left(X_{s}^{c}- X_{s}\right)-\int_{0}^{t}\left(X_{-}^{c}-X\right)\mathrm{d}Y\right| 
\\
& \leq & c\left(\left|Y_{0}\right|+\left|Y_{t}\right|\right)+\left|\sum_{0<s\leq t}\Delta Y_{s}\Delta\left(X_{s}^{c}- X_{s}\right) \right|+\left|\int_{0}^{t}\left(X_{-}^{c}-X\right)\mathrm{d}Y\right|. 
\end{eqnarray*}
Thus we get 
\begin{eqnarray*}
\lefteqn{ \sup_{0\leq t\leq T}\left|\int_{0}^{t}Y_{-}\mathrm{d}X^{c}-\int_{0}^{t}Y_{-}\mathrm{d}X-\left[X^{cont},Y^{cont}\right]_{t}\right|   } \\
&  \leq &  c\left(\left|Y_{0}\right|+\sup_{0\leq t\leq T}\left|Y_{t}\right|\right) + \sup_{0\leq t\leq T}\left|\sum_{0<s\leq t}\Delta Y_{s}\Delta\left(X_{s}^{c}- X_{s}\right) \right|
 + \sup_{0\leq t\leq T}\left|\int_{0}^{t}\left(X_{-}^{c}-X\right)\mathrm{d}Y\right|.
\end{eqnarray*}
Since $Y$ has càdlàg paths, it is bounded and hence $c\left(\left|Y_{0}\right|+\sup_{0\leq t\leq T}\left|Y_{t}\right|\right)\rightarrow0$
$\mathbb{P}$ a.s. as $c\downarrow0.$ 

Note that
\begin{eqnarray*}
\left|\Delta\left(X_{s}^{c}-X_{s}\right)\right| & = & \left|\left(X_{s}^{c}-X_{s}\right)-\left(X_{s-}^{c}-X_{s-}\right)\right|\leq 2c .\\
\end{eqnarray*}
Similarly, for $s\in\left[0;T\right],$ 
\[
\left|\Delta\left(X_{s}^{c}-X_{s}\right)\right|\leq\left|\Delta X_{s}^{c}\right|+\left|\Delta X_{s}\right|\leq\left(K_{T}+1\right)\left|\Delta X_{s}\right|.
\]
Thus we obtain that 
\begin{equation*}
\left|\Delta\left(X_{s}^{c}-X_{s}\right)\right|\leq\min\left\{ 2c,\left(K_{T}+1\right)\left|\Delta X_{s}\right|\right\} \leq\left(K_{T}+2\right)\min\left\{ c,\left|\Delta X_{s}\right|\right\} 
\end{equation*}
and using this, we estimate 
\begin{eqnarray*}
\sup_{0\leq t\leq T}\left|\sum_{0<s\leq t}\Delta Y_{s}\left(\Delta X_{s}^{c}-\mbox{\ensuremath{\Delta}}X_{s}\right)\right| & \leq & \sup_{0\leq t\leq T}\sqrt{\sum_{0<s\leq t}\left|\Delta\left(X_{s}^{c}-X_{s}\right)\right|^{2}}\sqrt{\sum_{0<s\leq t}\left|\Delta Y_{s}\right|^{2}} \\
& = & \sqrt{\sum_{0<s\leq T}\left|\Delta\left(X_{s}^{c}-X_{s}\right)\right|^{2}}\sqrt{\sum_{0<s\leq T}\left|\Delta Y_{s}\right|^{2}}\\
& \leq & \left(K_{T}+2\right)\sqrt{\sum_{0<s\leq T}\min\left\{ c^{2},\left|\Delta X_{s}\right|^{2}\right\} }\sqrt{\left[\Delta Y\right]_{T}}\rightarrow 0 \mbox{ } \mathbb{P} \mbox{ a.s.} \mbox{ as } c\downarrow 0.
\end{eqnarray*}
 In order to estimate 
\[
\sup_{0\leq t\leq T}\left|\int_{0}^{t}\left(X_{-}^{c\left(n\right)}-X_{-}\right)\mathrm{d}Y\right|
\]
let us decompose the semimartingale $Y$ into a local martingale $M$
and a locally finite variation process $A$ (note that they may depend
on the measure $\mathbb{P}),$ 
\[
Y=M+A.
\]
Let $\left(\tau\left(k\right)\right)_{k\geq1}$ be a sequence of stopping
times increasing to $+\infty$ such that $\left(M_{t\wedge\tau\left(k\right)}\right)_{t \geq 0}$
is a square integrable martingale. Using elementary estimate $\left(a+b\right)^{2}\leq2a^{2}+2b^{2}$
and Burkholder inequality, on the set $\Omega_{N}=\left\{ \omega\in\Omega:TV\left(A,\left[0;T\right]\right)\leq N\right\} $
we obtain 
\begin{eqnarray*}
\lefteqn{ \mathbb{E}\left[\sup_{0\leq t\leq T\wedge\tau\left(k\right)}\left|\int_{0}^{t}\left(X_{-}^{c}-X_{-}\right)\mathrm{d}Y\right|^{2};\Omega_{N}\right]}
\\
& \leq & 2\mathbb{E}\sup_{0\leq t\leq T\wedge\tau\left(k\right)}\left|\int_{0}^{t}\left(X_{-}^{c}-X_{-}\right)\mathrm{d}M\right|^{2}+2\left[\mathbb{E}\left|\int_{0}^{T}\left(X_{-}^{c}-X_{-}\right)\mathrm{d}A\right|^{2};\Omega_{N}\right]\\
& \leq  &
2\left(4c^{2}\mathbb{E}\left[M,M\right]_{T\wedge\tau\left(k\right)}+c^{2}N^{2}\right)=\left(8\mathbb{E}\left[M,M\right]_{T\wedge\tau\left(k\right)}+2N^{2}\right)c^{2}.\label{eq:est3}
\end{eqnarray*}
Let now $\left(c\left(n\right)\right)_{n\geq1}$ be such a sequence
that $\sum_{n=1}^{\infty}c\left(n\right)^{2}<+\infty.$ We have
\begin{eqnarray*}
\mathbb{E}\left[\sum_{n=1}^{\infty}\sup_{0\leq t\leq T\wedge\tau\left(k\right)}\left|\int_{0}^{t}\left(X_{-}^{c\left(n\right)}-X_{-}\right)\mathrm{d}Y\right|^{2};\Omega_{N}\right] & = & \sum_{n=1}^{\infty}\mathbb{E}\left[\sup_{0\leq t\leq T\wedge\tau\left(k\right)}\left|\int_{0}^{t}\left(X_{-}^{c\left(n\right)}-X_{-}\right)\mathrm{d}Y\right|^{2};\Omega_{N}\right]\\
 & \leq & \left(8\mathbb{E}\left[M,M\right]_{T\wedge\tau\left(k\right)}+2N^{2}\right)\sum_{n=1}^{\infty}c\left(n\right)^{2}\\
 & < & +\infty.
\end{eqnarray*}
Hence, the sequence $\sup_{0\leq t\leq T\wedge\tau\left(k\right)}\left|\int_{0}^{t}\left(X_{-}^{c\left(n\right)}-X_{-}\right)\mathrm{d}Y\right|^{2}, n=1,2,...,$
converges to $0$ on the set $\Omega_{N}.$ Since $\Omega=\bigcup_{N\geq1}\Omega_{N},$
we get that $\sup_{0\leq t\leq T\wedge\tau\left(k\right)}\left|\int_{0}^{t}\left(X_{-}^{c\left(n\right)}-X_{-}\right)\mathrm{d}Y\right|^{2}$
converges $\mathbb{P}$ a.s. to $0.$ Finally, since $\tau\left(k\right)\rightarrow+\infty$
we get that $\sup_{0\leq t\leq T}\left|\int_{0}^{t}\left(X_{-}^{c\left(n\right)}-X_{-}\right)\mathrm{d}Y\right|^{2}$
converges $\mathbb{P}$ a.s. to $0.$

\end{proof}

\section{$p-$variation of the sequence $\left(X^{c}\right)_{c>0}$}

Let us fix $p>0$ and$0\leq a<b.$ We will consider $p-$variation
of the càdlàg process $X$ defined (pathwise) in the three following
ways. 
\begin{enumerate}
\item 
\[
v_{p}^{\left(1\right)}\left(X,\left[a;b\right]\right)=\limsup_{\delta\rightarrow0}\sup_{n}\sup_{\substack{a\leq t_{0}<t_{1}<...<t_{n}\leq b\\
t_{i}-t_{i-1}\leq\delta\text{ for }i=1,2,...,n
}
}\sum_{i=1}^{n}\left|X_{t_{i}}-X_{t_{i-1}}\right|^{p},
\]

\item 
\[
v_{p}^{\left(2\right)}\left(X,\left[a;b\right]\right)=\limsup_{n\rightarrow\infty}\sup_{\substack{a\leq t_{0}<t_{1}<...<t_{n}\leq b\\
t_{i}-t_{i-1}\leq\delta_{n}\text{ for }i=1,2,...,n
}
}\sum_{i=1}^{n}\left|X_{t_{i}}-X_{t_{i-1}}\right|^{p},
\]
 where $\delta_{n}\rightarrow0$ as $n\rightarrow\infty,$ and 
\item 
\[
v_{p}^{\left(3\right)}\left(X,\left[a;b\right]\right)=\limsup_{n\rightarrow\infty}\sum_{i=1}^{n}\left|X_{t_{i}^{\left(n\right)}}-X_{t_{i-1}^{\left(n\right)}}\right|^{p},
\]
 where $\pi^{\left(n\right)}=\left\{ a\leq t_{0}^{\left(n\right)}<t_{1}^{\left(n\right)}<...<t_{n}^{\left(n\right)}\leq b\right\} $
is a nested sequence of partitions, $\pi^{\left(n\right)}\subset\pi^{\left(n+1\right)},$
$n=1,2,...,$ with $\max_{i=1,2,...,n}\left(t_{i}^{\left(n\right)}-t_{i-1}^{\left(n\right)}\right)\rightarrow0.$ 
\end{enumerate}
For example, for $W$ being a standard Wiener process, we have 
\[
v_{p}^{\left(1\right)}\left(W,\left[a;b\right]\right)=\left\{ \begin{array}{c}
\infty\text{ if }p\leq2,\\
0\text{ if }p>2
\end{array}\right.\text{a.s.}
\]
 (cf. \cite{Levy:1940}); 
\[
v_{p}^{\left(2\right)}\left(W,\left[a;b\right]\right)=\left\{ \begin{array}{c}
\infty\text{ if }p<2;\\
0\text{ if }p>2
\end{array}\right.\text{a.s.},
\]
 $v_{2}^{\left(2\right)}\left(W,\left[a;b\right]\right)$ is a.s.
finite but not fixed when $\delta_{n}=O\left(1/\ln\left(n\right)\right)$\ (cf.
\cite{delaVega:1974}) and $v_{2}^{\left(2\right)}\left(W,\left[a;b\right]\right)$
is a.s. finite and equal $b-a$ when $\delta_{n}=o\left(1/\ln\left(n\right)\right)$\ (cf.
\cite{Dudley:1973}). Finally, for $v^{(3)}$ one has (cf. \cite{Levy:1940})
\[
v_{p}^{\left(3\right)}\left(W,\left[a;b\right]\right)=\left\{ \begin{array}{c}
\infty\text{ if }p<2;\\
b-a\text{ if }p=2;\\
0\text{ if }p>2.
\end{array}\right.\text{a.s.}
\]

Let us now assume that for any $p>0,$ $p-$variation, $v_{p},$ is
defined by one of these three formulas or is defined as the limit
of equidistant $p-$variations, defined by formula (\ref{eq:realized}).
For processes $X$ and $X^{c}\in{\mathcal{X}}^{c}$ one may look at
the decomposition 
\[
X=X^{c}+\left(X-X^{c}\right)
\]
 as a decomposition into the sum of a process $X^{c}$ with locally
bounded total variation, uniformly approximating $X,$ and a ''noise''
$X-X^{c}$\ with small (smaller than $c$) amplitude but possibly
{}``wild'' behaviour with infinite total variation. Notice further
that since $X^{c}$ has finite total variation, it may be expressed
as the sum 
\[
X_{t}^{c}=\left(X^{c}\right)_{t}^{cont}+\sum_{0<s\leq t}\Delta X_{s}^{c},
\]
 where $\left(X^{c}\right)^{cont}$ is a continuous process. This
simple observation allows us to state 

\begin{thm} For any $p>1,$ the càdlàg process $X$ and a sequence
$\left(X^{c}\right)_{c>0}$ such that $X^{c}\in{\mathcal{X}}^{c}$
we have 
\begin{equation}
v_{p}\left(X^{c},\left[0;T\right]\right)=v_{p}\left(\sum_{0<s\leq\cdot}\Delta X_{s}^{c},\left[0;T\right]\right)=\sum_{0<s\leq T}\left|\Delta X_{s}^{c}\right|^{p}<+\infty\label{eq:pvar1}
\end{equation}
 and 
\[
\lim_{c\downarrow0}v_{p}\left(X^{c},\left[0;T\right]\right)=\lim_{c\downarrow0}\sum_{0<s\leq T}\left|\Delta X_{s}^{c}\right|^{p}=\sum_{0<s\leq T}\left|\Delta X_{s}\right|^{p}
\]
 which may be finite or infinite. Moreover, for any such $p\geq1$
that $v_{p}\left(X,\left[0;T\right]\right)=+\infty$ a.s., 
\begin{equation}
v_{p}\left(X-X^{c},\left[0;T\right]\right)=v_{p}\left(X,\left[0;T\right]\right)=+\infty\mbox{ a.s.}\label{eq:pvar2}
\end{equation}
 \end{thm} 

\begin{proof} Denote 
\[
Z^{c}=\left(X^{c}\right)^{cont}\mbox{ and \ensuremath{D^{c}=\sum_{0<s\leq\cdot}\left|\Delta X_{s}^{c}\right|.}}
\]
 To prove equality (\ref{eq:pvar1}) we start with the simple observation
that for $c>0\mbox{ and }p>1$ 
\[
v_{p}\left(Z^{c},\left[0;T\right]\right)=0.
\]
 This follows from continuity of $Z^{c}$ and finiteness of its total
variation. Indeed, for any partition 
\[
\pi=\left\{ 0\leq t_{0}<t_{1}<...<t_{n}\leq T\right\} 
\]
 we have 
\begin{eqnarray*}
\sum_{i=1}^{n}\left|Z_{t_{i}}^{c}-Z_{t_{i-1}}^{c}\right|^{p} & \leq & \max_{i=1,2,...,n}\left|Z_{t_{i}}^{c}-Z_{t_{i-1}}^{c}\right|^{p-1}\left(\sum_{i=1}^{n}\left|Z_{t_{i}}^{c}-Z_{t_{i-1}}^{c}\right|\right)\\
 & \leq & \omega\left(\mbox{mesh}\left(\pi\right),Z^{c}\right)^{p-1}TV\left(Z^{c},\left[0;T\right]\right),
\end{eqnarray*}
 where 
\[
\omega\left(h,f\right):=\sup_{\left|x-y\right|\leq h}\left|f\left(x\right)-f\left(y\right)\right|
\]
 denotes the modulus of continuity of a function $f$ and $\mbox{mesh}\left(\pi\right):=\max_{i=1,2,...,n}\left(t_{i}-t_{i-1}\right).$
Hence, by the definition of $p-$variation, the above inequality and
by continuity of $Z^{c},$ 
\begin{eqnarray*}
v_{p}\left(Z^{c},\left[0;t\right]\right) & \leq & v_{p}^{\left(1\right)}\left(Z^{c},\left[0;T\right]\right)\\
 & \leq & \limsup_{\delta\rightarrow0}\omega\left(\delta,Z^{c}\right)^{p-1}TV\left(Z^{c},\left[0;T\right]\right)=0.
\end{eqnarray*}
 Now observe that for any $p\geq1$ and any deterministic functions
$f$ and $g,$ 
\begin{equation}
v_{p}\left(f+g,\left[0;T\right]\right)^{1/p}\leq v_{p}\left(f,\left[0;T\right]\right)^{1/p}+v_{p}\left(g,\left[0;T\right]\right)^{1/p},\label{eq:mink}
\end{equation}
 which is a direct consequence of the Minkowski inequality. Using
(\ref{eq:mink}), we obtain 
\begin{eqnarray}
v_{p}\left(X^{c},\left[0;T\right]\right)^{1/p} & = & v_{p}\left(Z^{c}+D^{c},\left[0;T\right]\right)^{1/p}\leq v_{p}\left(Z^{c},\left[0;T\right]\right)^{1/p}+v_{p}\left(D^{c},\left[0;T\right]\right)^{1/p}\nonumber \\
 & = & v_{p}\left(D^{c},\left[0;T\right]\right)^{1/p}\label{eq:mink1}
\end{eqnarray}
 and, similarly, 
\begin{eqnarray}
v_{p}\left(D^{c},\left[0;T\right]\right)^{1/p} & = & v_{p}\left(X^{c}-Z^{c},\left[0;T\right]\right)^{1/p}\leq v_{p}\left(X^{c},\left[0;T\right]\right)^{1/p}+v_{p}\left(-Z^{c},\left[0;T\right]\right)^{1/p}\nonumber \\
 & = & v_{p}\left(X^{c},\left[0;T\right]\right)^{1/p}.\label{eq:mink2}
\end{eqnarray}
 Equality (\ref{eq:pvar1}) follows now from (\ref{eq:mink1}) and
(\ref{eq:mink2}), and the simple observation that $v_{p}\left(D^{c},\left[0;T\right]\right)=\sum_{0<s\leq T}\left|\Delta X_{s}^{c}\right|^{p}<+\infty$. 

The limit 
\[
\lim_{c\downarrow0}v_{p}\left(X^{c},\left[0;T\right]\right)=\lim_{c\downarrow0}\sum_{0<s\leq T}\left|\Delta X_{s}^{c}\right|^{p}=\sum_{0<s\leq T}\left|\Delta X_{s}\right|^{p}
\]
 follows now from (\ref{eq:pvar1}) and uniform, dominated convergence
$X_{t}^{c}\rightrightarrows X_{t}$ as $c\downarrow0$ (recall that
for any $0\leq s\leq T,$ $\left|\Delta X_{s}^{c}\right|\leq K_{T}\left|\Delta X_{s}\right|;$
also note that the bound $\left|\Delta X^{c}\right|\leq\left|\Delta X\right|+2c$
and a.s. pointwise convergence $X_{t}^{c}\rightarrow X_{t}$ as $c\downarrow0$
are sufficient).

In order to prove (\ref{eq:pvar2}) one now may utilize finiteness
of $v_{p}\left(X^{c},\left[0;T\right]\right)$ and (\ref{eq:mink}).
\end{proof} 

We have just proved that the sequence of the limits of equidistant
$p-$variations, $v_{p}^{\left(0\right)}\left(X^{c},\left[0;T\right]\right)_{c>0},$
tends (pathwise) to $\sum_{0<s\leq T}\left|\Delta X_{s}\right|^{p}$.
It is interesting to compare this result with the results of \cite{Lepingle:1976}
and \cite{Jacod:2008} on the limits (in probability) of equidistant
$p-$variations, $\tilde{v}_{p}^{\left(0\right)}\left(X;\left[0;T\right]\right),$
defined as 
\begin{equation}
\tilde{v}_{p}^{\left(0\right)}\left(X;\left[0;T\right]\right):=\left(\mathbb{P}\right)\lim_{n\rightarrow\infty}\sum_{i=1}^{n}\left|X_{iT/n}-X_{\left(i-1\right)T/n}\right|,\label{eq:realized1}
\end{equation}
 when $X$ is a semimartingale on some filtered probability space
$\left(\Omega,{\mathcal{F}},F,\mathbb{P}\right).$ Since $X$ is a
semimartingale, it may be written as 
\begin{eqnarray*}
X_{t} & = & X_{0}+\int_{0}^{t}b_{s}\mathrm{d}s+\int_{0}^{t}\sigma_{s}\mathrm{d}W_{s}+\int_{0}^{t}\int_{E}\kappa\circ\delta\left(s,x\right)\left(\mu-\nu\right)\left(\mathrm{d}s,\mathrm{d}x\right)\\
 &  & +\int_{0}^{t}\int_{E}\kappa'\circ\delta\left(s,x\right)\mu\left(\mathrm{d}s,\mathrm{d}x\right),
\end{eqnarray*}
 where $W$ and $\mu$ are a Wiener process and a Poisson random measure
on $[0;+\infty)\times E$ with $(E,{\mathcal{E}})$ an auxiliary measurable
space on the space $\left(\Omega,{\mathcal{F}},F,\mathbb{P}\right)$
(for details cf. \cite[description of formula (1)]{AitSahaliaJacod:2009}.
The mentioned results state that (cf. \cite[description of formula (11)]{AitSahaliaJacod:2009})
\[
\tilde{v}_{p}^{\left(0\right)}\left(X;\left[0;T\right]\right)=\begin{cases}
\sum_{0<s\leq T}\left|\Delta X_{s}\right|^{p} & \mbox{for }p>2;\\
\int_{0}^{t}\sigma_{s}^{2}\mathrm{d}s+\sum_{0<s\leq T}\left|\Delta X_{s}\right|^{2} & \mbox{for }p=2;
\end{cases}
\]
 and for $p\in\left(0;2\right),$ 
\[
\left(\mathbb{P}\right)\lim_{n\rightarrow\infty}\left(\frac{T}{n}\right)^{1-p/2}\sum_{i=1}^{n}\left|X_{iT/n}-X_{\left(i-1\right)T/n}\right|=\pi^{-1/2}2^{p/2}\Gamma\left(\frac{p+1}{2}\right)\int_{0}^{t}\left|\sigma_{s}\right|^{p}\mathrm{d}s.
\]
 Thus, we observe that for $p\in\left(1;2\right]$ we obtain different
limits for $X$ than for $X^{c}.$ Otherwise as for the process with
locally finite total variation, the influence of small jumps and continuous
part of $X$ on $p-$variation can not be neglected for $p\leq2.$

\end{document}